\documentclass{article} 
\usepackage{mathmag}

\usepackage{amsmath,amsthm}     
\usepackage{graphicx}     
\usepackage{hyperref} 
\usepackage{url}
\usepackage{amsfonts} 

\usepackage{multicol}

\theoremstyle{theorem}
\newtheorem{theorem}{Theorem}

\theoremstyle{definition}

\newtheorem*{remark}{Remark}

\allowdisplaybreaks

\makeatletter
\@addtoreset{footnote}{page}
\makeatother

\begin{document}

\title{Finding Fibonacci in the Hyperbolic Plane}

\author{MurphyKate Montee\\               
\scriptsize Carleton College\\    
Northfield, MN 55057 \\
mmontee@carleton.edu}                      

\maketitle

This paper is a result of a mathematical mistake I made seven years ago. As a graduate student, a math professor of mine asked me to compute the perimeters and areas of some hyperbolic polygons as an introduction to hyperbolic geometry. I sat down, drew some pictures, found a pattern, and gave him an answer. We were both happy, and went on with our study of hyperbolic spaces.

Last year, while I was teaching hyperbolic geometry for the first time, I remembered that problem. I'd had fun solving it and it had a nice and neat solution, so I put it on the homework set for my students. They were much more careful than either my former self or my former professor, and they discovered what neither of us had: there was a mistake in my solution from seven years ago! The correct solution ended up being more of a challenge to prove than I had planned, but it was also a much more surprising and beautiful result.

\section{The Problem}

Consider the family of complexes $X_n$ made by gluing together equilateral triangles along their edges so that there are $n$ triangles at every vertex. If $n$ is $3, 4,$ or $5$, we get a finite complex; in fact, we get the three Platonic solids with triangular faces. If $n=6$, we get a copy of the Euclidean plane (tiled by equilateral triangles). 

\begin{figure}[h]
    \centering
    \includegraphics[width=.9\textwidth]{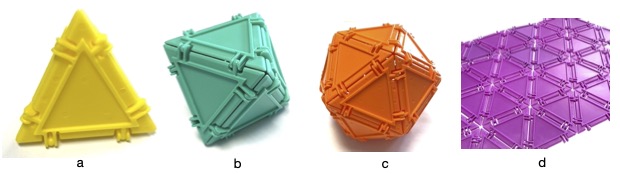}
    \caption{The complexes (a) $X_3,$ (b) $X_4,$ (c) $X_5$, (d) part of $X_6$. Note that for $n\geq 6,$, $X_n$ is infinite. All models made from Geometiles\textregistered.}
    \label{fig:family_of_spaces}
\end{figure}

\begin{figure}
    \centering
    \includegraphics[width=.4\textwidth]{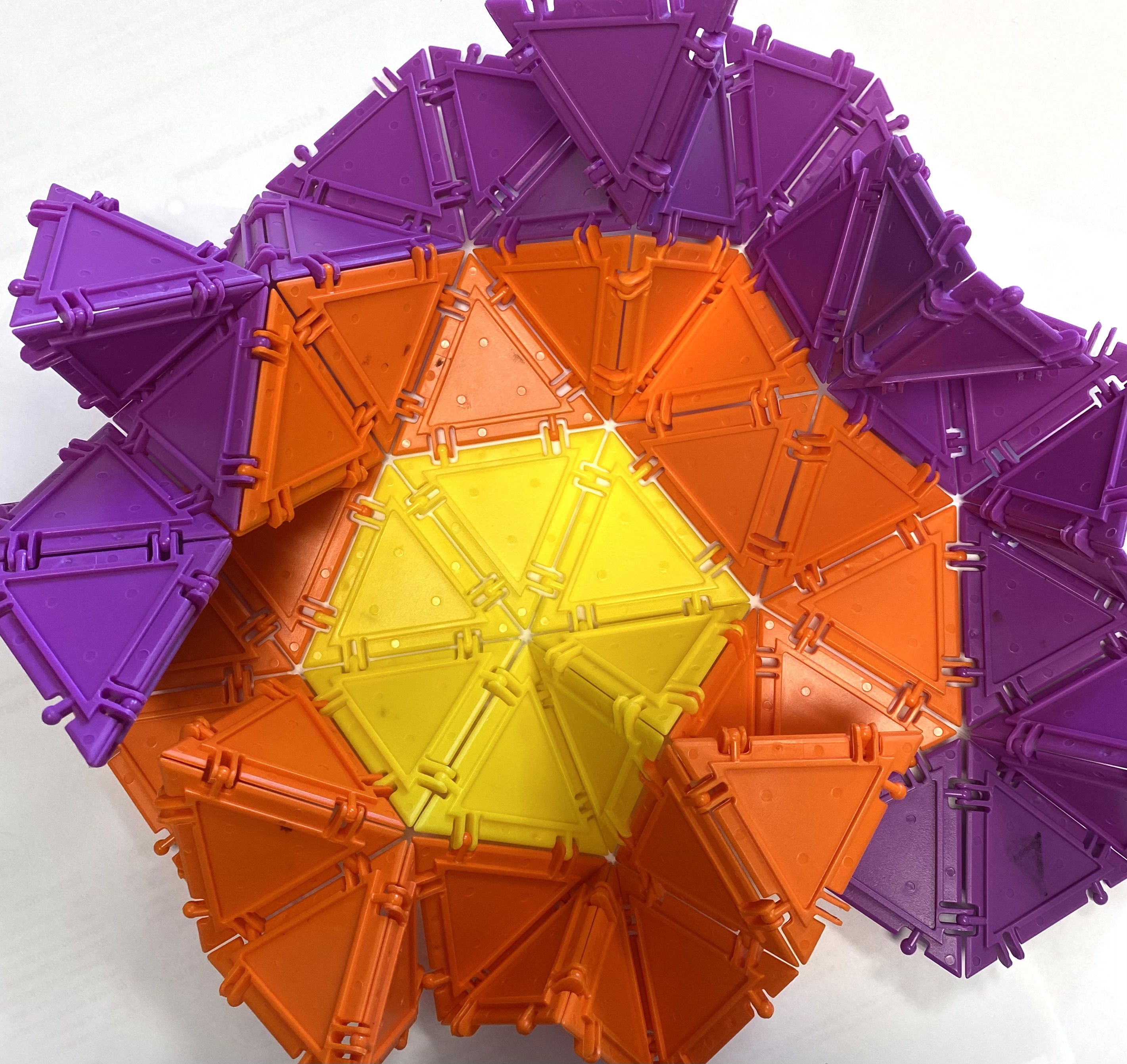}
    \caption{Part of $X_7$, also made from Geometiles \textregistered.}
    \label{fig:X7}
\end{figure}

\begin{figure}
    \centering
    \includegraphics[width=.5\textwidth]{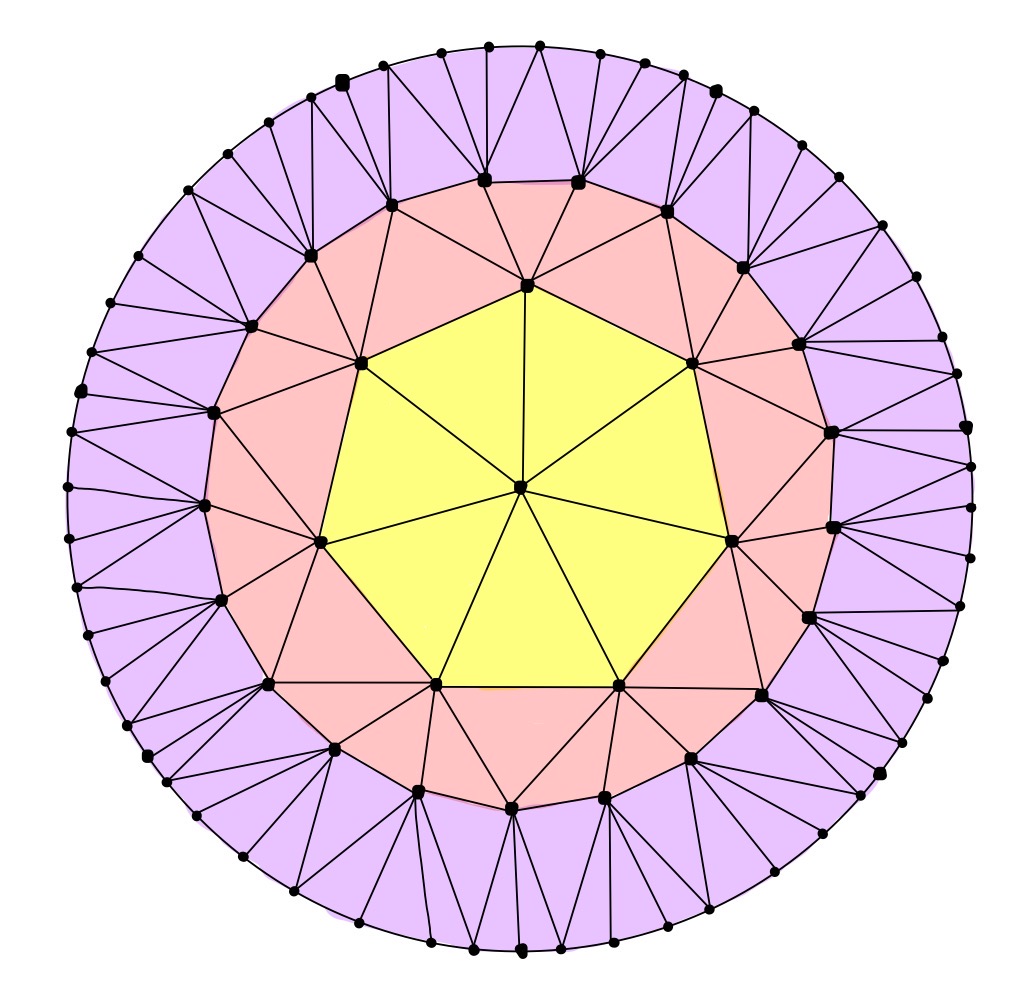}
    \caption{Part of $X_7$. (Not drawn to scale!) The yellow triangles lie in $D(1)$, the red in $D(2)$, and the purple in $D(3)$. The vertex $v$ is in the middle of the yellow region.}
    \label{fig:circles}
\end{figure}

If $n>6$, we get a strange and infinite surface that can't be drawn to scale in the Euclidean plane, though we can make it out of physical triangles. The result is a floppy and frilly object that refuses to lie flat, like a leaf of kale (see Figure \ref{fig:X7}). We can also draw it in the plane, if we accept a drawing that isn't accurate with respect to distances or angles (see Figure \ref{fig:circles}).

Let's scale our triangles so that each edge has length $1$. Pick any vertex, $v$. We can define a \textit{combinatorial disk in $X_n$ of radius $r$,} denoted $D_n(r)$, iteratively as follows: Pick any vertex $v$. Let $D_n(0) = \{v\}$, and let $D_n(r+1)$ be the union of $D_n(r)$ and all of the triangles that share a vertex with $D_n(r)$.

Let $P_n(r)$ denote the perimeter of the combinatorial disk in $X_n$. If we choose $n=6$, we can see that $D_6(r)$ is a hexagon with edge length $r$, so $P_6(r) = 6r$.  

When $n=7$, it's a more complicated story. Figure \ref{fig:circles} shows $D_7(1), D_7(2), D_7(3)$. Counting the perimeter of each combinatorial disk, you get the following sequence:

    \[
        \begin{tabular}{r||c|c|c|c|c|c}
        $r$           &1 &2  &3 &4 &5 &6\\
        \hline
        $P_7(r)$        &7 &21 &56 &91 &385 &938
        \end{tabular}
    \]

Do you see the pattern? It's not immediately obvious! But if we divide every entry by 7, we get the numbers $1, 3, 8, 21, 55, 144,$ and now the pattern emerges: these are every other Fibonacci number.

\begin{theorem}\label{thm:fib_numbers}
Let $F_n$ denote the $n$-th Fibonacci number. Then $P_7(r) = 7F_{2r}.$
\end{theorem}

The goal of the next section is to prove this theorem. 

\begin{remark}
To ease notation, for the rest of the paper we will write $D(r)$ to denote $D_7(r)$ and $P(r)$ to denote $P_7(r)$.
\end{remark}

\section{Fibonacci Appears}

We will prove Theorem \ref{thm:fib_numbers} by induction on the radius, $r$. We can establish the base case by checking $D(1)$. For the inductive step, we'll count the number of vertices in $\partial D(r)$ instead of the number of edges. 

\begin{figure}
    \centering
    \includegraphics[width=.7\textwidth]{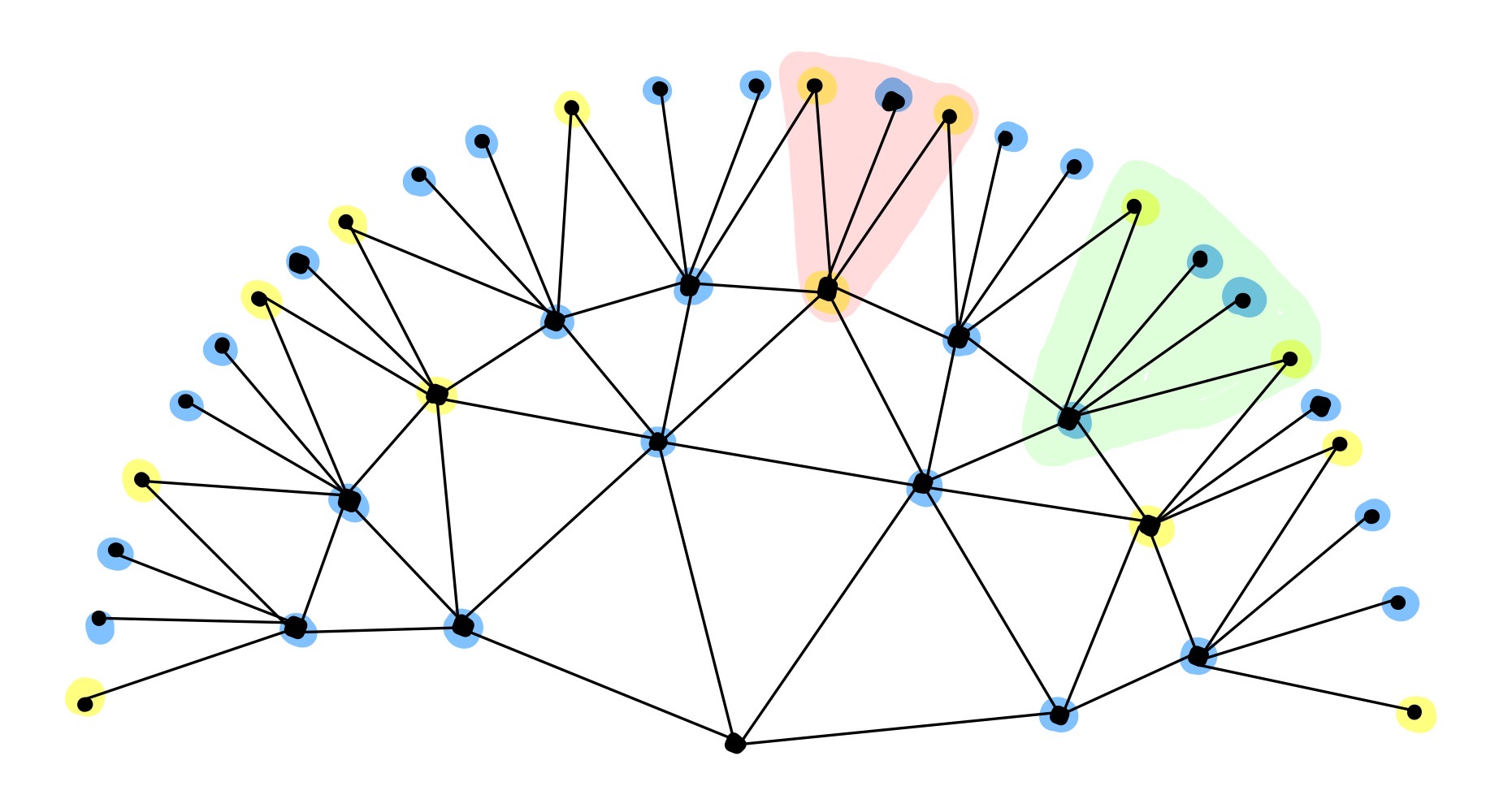}
    \caption{Part of $X_7.$ Vertices are colored blue and yellow as in the proof. The red and green regions each contain one vertex in $\partial D(2)$ and the vertices in $\partial D(3)$ that belong or partially belong to it.}
    \label{fig:inductive_step}
\end{figure}

Consider Figure \ref{fig:inductive_step}. Every vertex in $\partial D(r)$ is adjacent to 7 triangles in total. Either 2 or 3 of those triangles lie in $D(r)$. We'll partition the vertices in $\partial D(r)$ into two sets; the \textit{blue} vertices are those that are adjacent to 2 triangles in $D(r)$, and the \textit{yellow} vertices are those that are adjacent to 3 triangles in $D(r)$. 
Let $B_r$ be the number of blue vertices in $\partial D(r)$ and let $Y_r$ be the number of yellow vertices. 

Two vertices are \textit{adjacent} if they are separated by an edge. Every vertex in $\partial D(r+1)$ is adjacent to exactly 1 or 2 vertices in $\partial D(r)$. We'll say that a vertex $v$ in $\partial D(r+1)$ \textit{belongs to} a vertex $v' \in \partial D(r)$ if $v'$ is the only vertex in $D(r)$ adjacent to $v$. We'll say $v$ \textit{partially belongs to} $v'$ if $v'$ is one of the two vertices in $D(r)$ adjacent to $v$. 

Then the total number of vertices in $D(r+1)$ is the sum of 
    \[
        \#\{\mbox{vertices that belong to }v\} + \frac{1}{2}\#\{\mbox{vertices that partially belong to }v\},
    \]
taken over every vertex in $\partial D(r)$.

For every blue vertex $b$ in $\partial D(r)$, there are 2 blue vertices in $\partial D(r+1)$ that belong to $b$ and 2 yellow vertices that partially belong to $b$. For every yellow vertex $y$ in $\partial D(r),$ $\partial D(r+1)$ contains 1 blue vertex that belongs to $y$ and 2 yellow vertices that partially belong to $y$. In summary, we have that $P(r+1) = B_{r+1} + Y_{r+1}$, where:
    \[
        B_{r+1} = 2B_r + 1Y_r \mbox{ and } Y_{r+1} = 1B_r + 1Y_r.
    \]
Written in matrix form, that is:
     \[
    \begin{bmatrix} B_{r+1}\\Y_{r+1}\end{bmatrix}
        = \begin{bmatrix} 2 &1\\1&1\end{bmatrix}\begin{bmatrix} B_r\\Y_r\end{bmatrix}.
    \]
In particular, since $B_1 = 7$ and $Y_1 = 0$, we get 
    \begin{equation}\label{eqn:matrix_form}
    \begin{bmatrix} B_{r+1}\\Y_{r+1}\end{bmatrix}
        = \begin{bmatrix} 2 &1\\1&1\end{bmatrix}^r\begin{bmatrix} 7\\0\end{bmatrix} = 7\cdot \begin{bmatrix} 2 &1\\1&1\end{bmatrix}^r\begin{bmatrix} 1\\0\end{bmatrix}.
    \end{equation}
    
On the other hand, the Fibonacci numbers also satisfy a linear recurrence relation, namely: 
    \[
    \begin{bmatrix} F_{r+1}\\F_{r}\end{bmatrix}
        = \begin{bmatrix} 1 &1\\1&0\end{bmatrix}^r\begin{bmatrix} 1\\0\end{bmatrix}.
    \]
Since we're only interested in even order Fibonacci numbers, we have
        \[
    \begin{bmatrix} F_{2r+1}\\F_{2r}\end{bmatrix}
        = \begin{bmatrix} 1 &1\\1&0\end{bmatrix}^{2r}\begin{bmatrix} 1\\0\end{bmatrix}
        = \begin{bmatrix} 2 &1\\1&1\end{bmatrix}^r\begin{bmatrix} 1\\0\end{bmatrix}.
    \]
This differs from the result of Equation \ref{eqn:matrix_form} by a factor of 7, so we get $7F_{2r+1} = B_{r+1}$ and $7F_{2r} = Y_{r+1}$. Therefore: 
    \[
        7F_{2r+2} = 7(F_{2r+1}+F_{2r}) = B_{r+1} + Y_{r+1} = P(r+1).
    \]

\section{Hyperbolicity, Exponential Growth, and Isoperimetry}
The title of this paper promised that we'd find Fibonacci in the hyperbolic plane, so where is the hyperbolic plane in all this?

One of the standard models of the hyperbolic plane is the \textit{Poincaré disk}. In this model, the `plane' is the interior of the unit disk, and geodesics are arcs of Euclidean circles or lines that cross the boundary of the disk at right angles.

Many results in hyperbolic geometry are extremely counterintuitive at first. One of the most surprising facts about hyperbolic triangles is that the sum of the interior angles is always strictly less than $\pi$ (in radians). Furthermore, the hyperbolic area of a hyperbolic triangle is determined entirely by the sum of the interior angles: A triangle $\Delta ABC$ with interior angles $\alpha, \beta, \gamma$ has area $\pi - (\alpha + \beta + \gamma)$. Two more surprising facts, which are consequences of this, are that (1) any two hyperbolic triangles which are similar must also be congruent, and (2) there exists a largest triangle, namely the so-called \textit{ideal triangle} whose angles are all 0. A good resource if you'd like to read more about hyperbolic geometry is \cite{stahl}.

\begin{figure}[h]
    \centering
    \includegraphics[width=.4\textwidth]{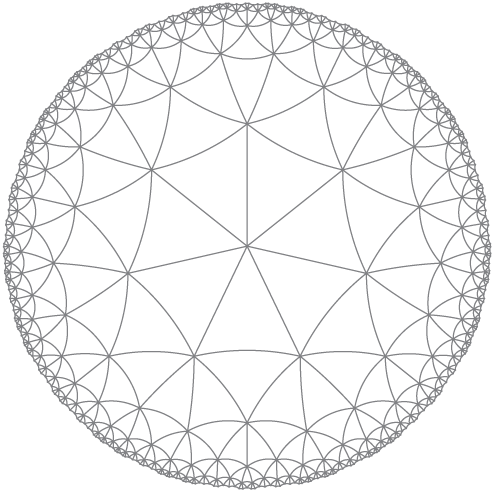}
    \caption{A tiling of the hyperbolic plane by equilateral triangles, with 7 triangles at every vertex. }
    \label{fig:7tiling}
\end{figure}

As we scale an equilateral hyperbolic triangle, the sum of the interior angles changes: the larger the triangle, the smaller the sum. Since the triangle is equilateral, it is also equiangular, so the angle at any one corner changes as well. Therefore by scaling our triangles to the appropriate size, we can find equilateral triangles with interior angles $2\pi/n$ for any $n>6$.  In particular, if we pick $n=7$ we will get equilateral triangles with interior angles perfectly sized to fit 7 triangles around a common vertex. These triangles have area equal to $\pi - \frac{6\pi}{7} = \frac{\pi}{7}$, by the fact above. The perimeter of the triangle is a little more difficult to calculate. If we continue to place these equilateral triangles, with 7 at each vertex, we get a tiling of the hyperbolic plane, as illustrated in Figure \ref{fig:7tiling}.

We built $X_7$ out of Euclidean triangles with side length $1$, but we can interpret \ref{thm:fib_numbers} as a statement about the number of edges around the boundary of $D(r),$ rather than its perimeter. In fact, we can call the number of edges around the boundary of $D(r)$ the \emph{combinatorial perimeter} of $D(r)$, and we'll denote it $P(r)$. (Note that $P(r)$ will always be a constant multiple of the actual perimeter of $D(r)$). In this interpretation, it doesn't matter if we use Euclidean or hyperbolic triangles to build $X_7$. If we build $X_7$ out of hyperbolic triangles our combinatorial disks $D(r)$ approximate hyperbolic disks.

This implies (at least asymptotically) another counterintuitive fact about hyperbolic space: The perimeter of a hyperbolic circle grows exponentially with the radius. Compare this to the Euclidean case, where the perimeter of a circle grows linearly with its radius!

Another fun fact about the hyperbolic plane is that polygons (and more general regions) satisfy a \textit{linear isoperimetric inequality}; that is, if $R$ is a hyperbolic polygon with perimeter $|\partial R|$ and area $|R|$, then there is a constant $c$ (which does not depend on the choice of polygon) so that $|\partial{R}| \geq c |R|$. Again, let's compare this to the Euclidean case. Consider any Euclidean polygon $R$, and let $R_k$ denote a copy of $R$ scaled by $k$. Then $|R_k|$ is proportional to $k^2$, while $|\partial R_k|$ is proportional to $k$. Thus $\lim_{k\to \infty} |\partial R_k| / |R_k| = 0, $ so there does not exist a constant $c$ such that $|\partial R_k| / |R_k| \geq c$ for all $k$. (And we didn't even have to consider non-similar polygons!)

For this paper, we'll use an approximation of area that's easier to calculate with, called the \emph{combinatorial area} and denoted $A(r)$, where $A(r)$ is the number of triangles required to make up our disk $D(r)$. Note that if we make $X_7$ out of hyperbolic triangles, then $\frac{\pi}{7}A(r) = |D(r)|$. If we make $A(r)$ out of Euclidean triangles, then $\frac{\sqrt{3}}{4}A(r) = |D(r)|$. In general, the combinatorial area $A(r)$ will always be a constant multiple of the actual area, $|D(r)|$, in the same way that the combinatorial perimeter is a constant multiple of the actual perimeter. In particular, our disks will satisfy a linear isoperimetric inequality with respect to combinatorial area and perimeter if and only if they satisfy one (possibly with a different coefficient) with respect to actual area and perimeter.

Can we see this linear isoperimetric inequality in our combinatorial model?

\section{Fibonacci Strikes Again}
Let's find a formula for $A(r)$. By counting triangles in Figure \ref{fig:circles}, we get the following data:
    \[
        \begin{tabular}{r||c|c|c|c|c|c}
        $r$           &1 &2  &3 &4 &5 &6\\
        \hline
        $A(r) - A(r-1)$        &7 &28 &77 &315 &847 &2240
        \end{tabular}
    \]
This time the pattern isn't as easy to spot. Here it is:

\begin{theorem}\label{thm:fib_area}
Let $A(r)$ denote the combinatorial area of $D(r)$. Then we have $A(r) = 7(4F_{2r-2} + 3F_{2r-3} - 2).$
\end{theorem}

\begin{figure}
    \centering
    \includegraphics[width=.6\textwidth]{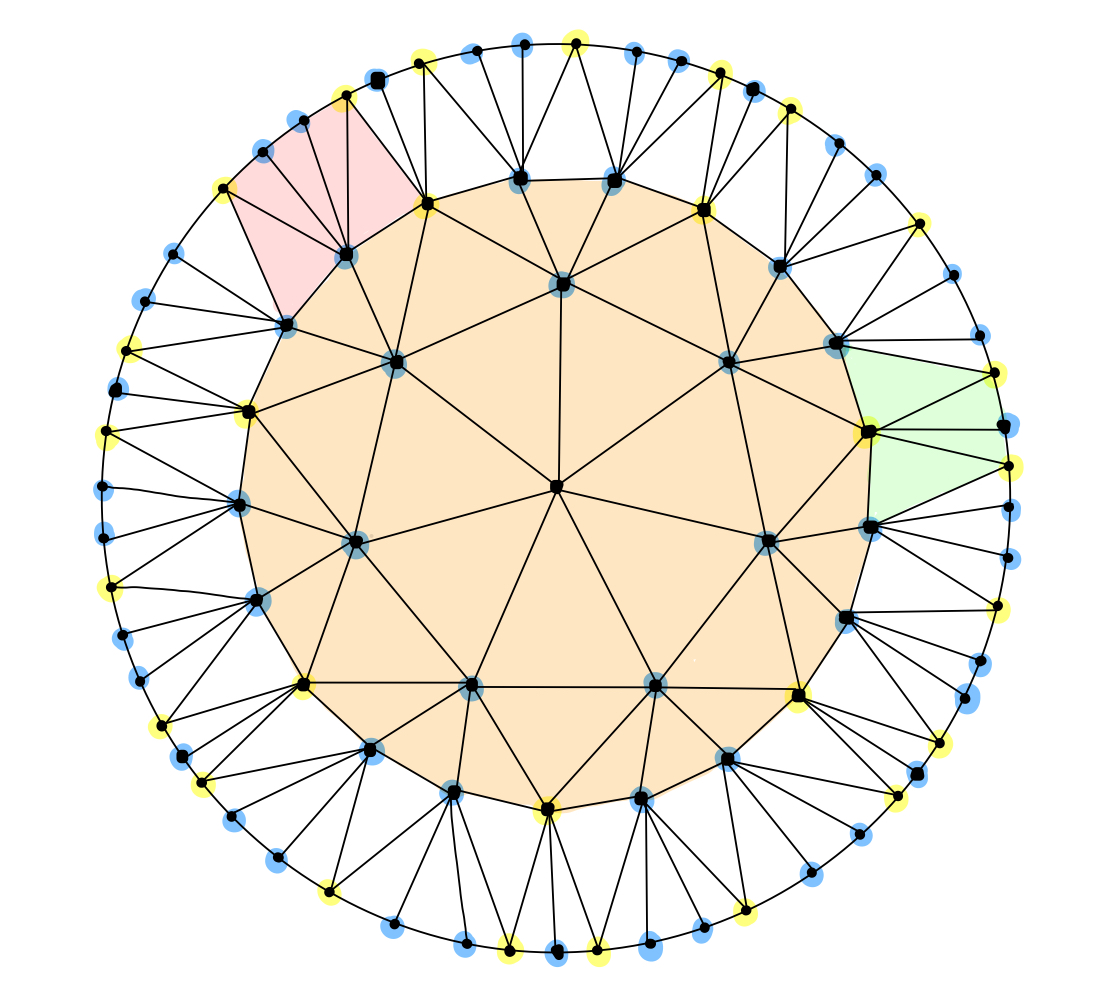}
    \caption{The orange region is $D(2).$ Blue vertices are adjacent to 5 triangles in $D(3)-D(2)$, as indicated in red. Yellow vertices are adjacent to 4, as indicated in green.}
    \label{fig:area}
\end{figure}

To prove the theorem, we'll enumerate triangles in $D(r+1)-D(r)$ using the blue and yellow vertices from before.  Every blue vertex in $\partial D(r)$ is adjacent to $5$ triangles in $D(r+1) - D(r)$. On the other hand, every yellow vertex in $\partial D(r)$ is adjacent to $4$ triangles in $D(r+1) - D(r)$. (See Figure \ref{fig:area}.) 

It's tempting to say that the change in area is therefore $5B_r + 4Y_r$, but be careful: some of the triangles in $D(r+1) - D(r)$ are adjacent to two vertices in $\partial D(r)$. In fact, every vertex is adjacent to exactly two of the double-counted triangles because they are precisely the triangles in $D(r+1) - D(r)$ that share an edge with a triangle in $D(r)$. But since the number of vertices in $\partial D(r)$ is $P_r = B_r + Y_r$, we get:
    \[
        A(r+1) - A(r) = 5B_r + 4Y_r -  P_r = 4B_r + 3Y_r
    \]
whenever $r\geq1$.

Since $B_r=7F_{2r-1}$ and $Y_r = 7F_{2r-2}$, this gives us
    \[
        A(r+1) - A(r) =7\left( 4F_{2r-1} + 3F_{2r-2} \right).\\
    \]
Therefore for $r\geq 1$, we have:
    \begin{align*}
        A(r) &= A(1) + \sum_{i=1}^{r-1}\left( A(i+1) - A(i)\right) \\
            &= 7 + 7\sum_{i=1}^{r-1} ( 4F_{2i-1} + 3F_{2i-2}).
    \end{align*}
This looks nasty, but since these are the Fibonacci numbers there is a lovely identity that simplifies the story. You can verify, either by induction or by using the recursion formula, that $\sum_{i=1}^r F_{2i} = F_{2r+1} - 1,$ and $\sum_{i=1}^rF_{2i-1} = F_{2r}.$

Therefore, for $r>1$ we get:
    \begin{align*}
        A(r) &=7\left(1 + 4F_{2r-2} + 3(F_{2r-3} -1 )  \right)\\
            &= 7(4F_{2r-2} + 3F_{2r-3} - 2).
    \end{align*}

This isn't quite as nice as the formula for the perimeter, but it does allow us to compare the perimeter and the area. In fact, we can use facts about the Fibonacci sequence to calculate $\lim_{r\to \infty} \frac{A(r)}{P(r)}$ exactly. You can show, using the matrix form of the Fibonacci recursion, that $F_r = \frac{\phi^r - \psi^r}{2}$, where $\phi = \frac{1+\sqrt{5}}{2}$ (the Golden Ratio), and $\psi = \frac{1-\sqrt{5}}{2}$. Then you can use this to calculate that $\lim_{r\to\infty} \frac{F_{2r-2}}{F_{2r}} = \frac{1}{\phi^2},$ and $\lim_{r\to\infty} \frac{F_{2r-3}}{F_{2r}} = \frac{1}{\phi^3}$. Therefore
    \[
        \lim_{r\to\infty}\frac{A(r)}{P(r)} = 4\frac{F_{2r-2}}{F_{2r}} + 3 \frac{F_{2r-3}}{F_{2r}} - 2 \frac{-2}{F_{2r}} = \frac{4}{\phi^2} + \frac{3}{\phi^3} = \sqrt{5}.
    \]
This is certainly different from the Euclidean world!

To show that $X_7$ approximates the hyperbolic plane, we would like to show that $P(r)\geq c A(r)$ for some constant $c>0$. Equivalently, we want to show that $A(r)/P(r) \leq \frac{1}{c}$ for some $c>0$. And since we've been haunted by the number 7 throughout this entire paper, it would be nice -- for the sake of narrative symmetry if nothing else -- if 7 shows up again. Luckily for us, it does.

\begin{theorem}
For any radius $r$, $P(r) \geq \frac{1}{7}A(r).$
\end{theorem}

Let's see why. Using our formulas from Theorems \ref{thm:fib_numbers} and \ref{thm:fib_area}, we have 
    \begin{align*}
        \frac{A(r)}{P(r)} &= \frac{4 F_{2r-2} + 3 F_{2r-3} - 2}{F_{2r}} \\
            &= 4\frac{F_{2r-2}}{F_{2r}} + 3\frac{F_{2r-3}}{F_{2r}} - \frac{2}{F_{2r}}.
    \end{align*}
But since $F_r$ is an increasing sequence, we know that $F_q / F_r \leq 1$ whenever $q \leq r$. So we get:
    \[
        \frac{A(r)}{P(r)} \leq 4 + 3 - 0 = 7,
    \]
just as we wanted.

\acknowledgments{The author would like the thank the anonymous reviewer for comments that improved the exposition of this paper.}


\begin{abstract}
    We use a combinatorial approximation of the hyperbolic plane to investigate properties of hyperbolic geometry such as exponential growth of perimeter and area of disks, and the linear isoperimetric inequality. This calculations give a surprising link to Fibonacci numbers.
\end{abstract}

MurphyKate Montee is an assistant professor of mathematics at Carleton College. Her research lies in geometric group theory, an area of math that uses geometric objects (like the combinatorial hyperbolic plane in this article) to understand groups.

\end{document}